\documentclass[12pt, reqno]{amsart}
\usepackage{amsmath, amstext, amsbsy, amssymb, amscd}

\newtheorem{theorem}{Theorem}[section]
\newtheorem{lemma}[theorem]{Lemma}

\newtheorem{proposition}[theorem]{Proposition}
\newtheorem{corollary}[theorem]{Corollary}
\theoremstyle{definition}

\theoremstyle{remark}
\newtheorem{remark}[theorem]{Remark}

\numberwithin{equation}{section}

\newcommand{\C}{ \mathbb C }
\newcommand{\ch}{{\rm ch}}
\newcommand{\Ext}{ \text{Ext} }

\newcommand{\la}{\lambda}

\newcommand{\osp}{\text{osp}}

\newcommand{\Z}{ \mathbb Z }
\newcommand{\e}{\chi}

\newcommand{\fg}{{\mathfrak g}}
\newcommand{\fb}{{\mathfrak b}}
\newcommand{\fh}{{\mathfrak h}}
\newcommand{\cU}{{\mathcal U}}
\newcommand{\cO}{{\mathcal O}}
\newcommand{\cF}{{\mathcal F}}

\begin{document}
\title[Representation theory of $\osp(2|2n)$]{A Fock space approach to representation theory of $\osp(2|2n)$}
\author[Shun-Jen Cheng]{Shun-Jen Cheng}
\address{Department of Mathematics, National
Taiwan University, Taipei, Taiwan
106}\email{chengsj@math.ntu.edu.tw}
\thanks{Partially supported by NSC of R.O.C.}

\author[Weiqiang Wang]{Weiqiang Wang}
\address{Department of Mathematics, University of Virginia,
Charlottesville, VA 22904} \email{ww9c@virginia.edu}
\thanks{Partially supported by NSA and NSF}

\author[R.B.~Zhang]{R.B.~Zhang}
\address{School of Mathematics and Statistics, University of Sydney, New South Wales
2006, Australia}\email{rzhang@maths.usyd.edu.au}
\thanks{Partially supported by Australian Research Council}

\maketitle
\date{}
%\tableofcontents
%
%
\begin{abstract}
A Fock space is introduced that admits an action of a quantum
group of type $A$ supplemented with some extra operators. The
canonical and dual canonical basis of the Fock space are computed
and then used to derive the finite-dimensional tilting and
irreducible characters for the Lie superalgebra $\osp(2|2n)$. We
also determine all the composition factors of the symmetric
tensors of the natural $\osp(2|2n)$-module.
\end{abstract}

%\tableofcontents

%%
%%
%%
%%
\section*{Introduction}

Supersymmetry usually manifests itself as concrete representations
of the relevant Lie superalgebras, and thus the representation
theory of Lie superalgebras plays an essential role in the study
of supersymmetry. A crucial difference between simple Lie algebras
and simple Lie superalgebras is that the categories of finite
dimensional representations of the latter are in general not
semi-simple \cite{K1, K2}. Partially inspired by earlier works of
Lascoux-Leclerc-Thibon \cite{LLT} and Serganova \cite{Se}, Brundan
in \cite{B1} developed a new approach to the representation theory
of the general linear superalgebra (and in \cite{B3} for
$\mathfrak q(n)$), where he formulated the Kazhdan-Lusztig theory
in terms of canonical bases of a Fock space. The Fock space
approach enabled Brundan to establish, among other results, a
conjecture of \cite{VZ} in the affirmative. The Fock space
approach is not only technically powerful, it also leads to a new
conceptual framework relating the representations of general
linear algebras to those of general linear superalgebras via a
fundamental duality \cite{CWZ}.

The main purpose of this paper is to extend Brundan's Fock space
approach of Kazhdan-Lusztig theory to the Lie superalgebras
$\osp(2|2n)$. After reviewing some background materials in
Section~\ref{sec:pre}, we introduce in Section \ref{Fock-space} a
Fock space with the action of a quantum group of type $A$
%which is the quantized universal enveloping algebra of the
%(one-sided) infinite dimensional general linear Lie algebra
supplemented with some extra generators. We remark that the action
of the quantum group arises from the two opposite
comultiplications. In contrast to \cite{B1}, a weight space of the
Fock space here may correspond to a unique block or sometimes to
two blocks. A bar involution is defined and the canonical and dual
canonical bases of the Fock space are worked out explicitly.
Theorem~\ref{main} shows that the canonical basis elements when
specialized to $q=1$ correspond to the finite dimensional tilting
modules of $\osp(2|2n)$ while the standard monomials correspond to
the Kac modules. This together with the general theory of tilting
modules \cite{So, B2} gives an explicit determination of the
composition factors of Kac modules (Corollary~\ref{vanderjeugt}).

We point out that the finite dimensional irreducible
representations of $\osp(2|2n)$ were studied in the work of van
der Jeugt \cite{V}, who in particular established the composition
factors of Kac modules in a completely different way. This readily
implies a Bernstein-Leites type character formula for such
representations as noted in \cite{V}. Also the projective covers
in the category of finite-dimensional $\osp(2|2n)$-modules were
understood in \cite{Zou}. The new Fock space approach here appears
to be conceptually interesting and simple, and it is our hope that
it may provide new insights into the representation theory of
other Lie superalgebras.

In spite of the power of the Kazhdan-Lusztig-Brundan theory, the
structure of various naturally constructed modules of Lie
superalgebras often remains unclear. The skew-symmetric tensor of
the natural $\osp(2|2n)$-module is easily seen to be irreducible.
In contrast to the classical Lie algebra setup, the symmetric
tensors of the natural module of $\osp(2|2n)$ is not completely
reducible in general and it is a rather nontrivial problem to
determine the composition factors. In Section~\ref{sec:tensor} we
offer a complete solution to this problem. There exists a
surjective homomorphism via a Laplacian operator $\Delta$ from the
$k$-th symmetric tensor to the $(k-2)$-th symmetric tensor for
each $k$. We show that the kernel of $\Delta$ has $1, 3,$ or $2$
composition factors depending on whether $k \le n$, $n <k\le 2n$,
or $k>2n$. The simplest case when $k\le n$ can be also found in
\cite{Lee}.

 \vspace{.1in}
\noindent {\bf Acknowledgments.} We thank all three host
institutions of the authors and NCTS-Taipei office for the
hospitality and support.

\section{Preliminaries} \label{sec:pre}
In this section we present some background material for the use in
later sections.
\subsection{Lie superalgebra $\osp(2|2n)$}
Throughout this paper, we shall denote by $\mathfrak g$ the Lie
superalgebra $\osp(2|2n)$ whose standard Dynkin diagram together
with the simple roots is given by:
\begin{equation*}%standard-dynkin%
%\vspace*{-8ex}$
%\begin{array}{c c}
\setlength{\unitlength}{0.16in}
\begin{picture}(20,1.5)
\put(0,1){\makebox(0,0)[c]{$\bigotimes$}}
\put(3.4,1){\makebox(0,0)[c]{$\bigcirc$}}
\put(10.2,1){\makebox(0,0)[c]{$\bigcirc$}}
\put(13.25,1){\makebox(0,0)[c]{$\bigcirc$}}
\put(0.5,1){\line(1,0){2.4}} \put(3.82,1){\line(1,0){2}}
\put(6.5,1){$\ldots$} \put(8.2,1){\line(1,0){1.5}}
\put(10.6,0.75){$\Longleftarrow$}\put(12.0,0.75){$=$}
\put(0,0){\makebox(0,0)[c]{$\epsilon-\delta_1$}}
\put(3.4,0){\makebox(0,0)[c]{$\delta_1-\delta_2$}}
\put(10,0){\makebox(0,0)[c]{$\delta_{n-1}-\delta_n$}}
\put(13.5,0){\makebox(0,0)[c]{$2\delta_{n}$}}
\end{picture}
\vspace{.2cm}
\end{equation*}

\noindent Here $\epsilon-\delta_1$ is odd. The set of positive
roots is a union of the even and the odd ones: $\Delta^+
=\Delta^+_0 \cup \Delta^+_1$. Denote by $\{e_i, f_i, h_i\}$ for
$i=0, 1, \dots, n$ the corresponding Chevalley generators of
$\mathfrak g$. Let $\rho=-n\epsilon + \sum_{i=1}^n
(n-i+1)\delta_i$, which is half the graded sum of the positive
roots of ${\mathfrak g}$. Let ${\fb}$ be the Borel subalgebra, and
${\fh}\subset \fb$ be the Cartan subalgebra of ${\mathfrak g}$
compatible with the above choice of the simple roots. The space
$\mathfrak h^*$ is endowed with a non-degenerate symmetric
bilinear form
\begin{eqnarray*} (\epsilon, \, \epsilon)=1, & (\delta_i, \,
\delta_j) = -\delta_{i j}, & (\epsilon, \, \delta_i)=0, \quad
\forall i, j.
\end{eqnarray*}
A weight $\mu$ is called {\em atypical} if there is an odd
positive root $\gamma=\epsilon-\delta_i$ or $\epsilon+\delta_i$
for some $i$ such that $(\mu+\rho, \gamma)=0$, and is called {\em
typical} otherwise (cf. \cite{K2}).

Denote by $\cO^+$ the category of finite dimensional $\Z_2$-graded
${\mathfrak g}$-modules of integral weights, that is, the weights
of every module belong to the $\Z$-span of $\epsilon$, $
\delta_i$, $i\ge 1$.  We denote by $X_+^{1|n}$ the set of the
dominant integral weights of ${\mathfrak g}$, namely,
\begin{eqnarray*}
X_+^{1|n}= \left.\left\{ \lambda=\lambda_{-1}\epsilon +
\sum_{i=1}^n\lambda_i\delta_i \right| \lambda_i\in \Z, \forall i;\
\, \lambda_1\ge \lambda_2\ge \dots\ge \lambda_n\ge 0 \right\}.
\end{eqnarray*}
The Lie superalgebra $\fg$ has a $\Z$-grading $\fg =\fg_- +\fg_0
+\fg_+$, where $\fg_0$ is its even subalgebra and $\fg_+$
(respectively $\fg_-$) is the subalgebra spanned by the odd
positive (respectively negative) root vectors. The Kac module is
defined as $K(\la) : = U(\fg) \otimes_{U(\fg_0 +\fg_+)} L^0(\la)$,
where $L^0(\la)$ denotes the irreducible $\fg_0$-module of highest
weight $\la$ (extended trivially to $\fg_0 +\fg_+$). Then the
irreducible module $L(\lambda)$ and Kac module $K(\lambda)$ with
highest weight $\lambda$ belong to $\cO^+$ if and only if
$\lambda\in X_+^{1|n}$. It is known that $K(\la) =L(\la)$ if $\la$
is typical.
\subsection{The Bruhat order} \label{bruhat}

The Weyl group $W$ of ${\mathfrak g}$ is defined to be the Weyl
group of the $\text{sp}(2n)$ subalgebra, which is generated by the
reflections corresponding to the even simple roots of $\fg$. If an
element $\lambda\in X_+^{1|n}$ is atypical with respect to an odd
positive root $\gamma$, we define
\begin{eqnarray}
\lambda^{\texttt L}:=w(\lambda+\rho -k \gamma) -\rho.
\label{L-operator}
\end{eqnarray}
Here $k$ is the smallest positive integer such that $\lambda+\rho
- k \gamma$ is $\text{sp}(2n)$-regular in the sense that
$(\lambda+\rho - k \gamma, \alpha_i)\ne 0$, $\forall i\ge 1$, and
$w$ is the unique element in the Weyl group of ${\mathfrak g}$
rendering $\lambda^{\texttt L}$ dominant. For example, if
$\la=\epsilon$, then $\lambda^{\texttt
L}=-\epsilon+\delta_1+\delta_2$. A more explicit description of
this ``${\texttt L}$-operator" will be given below. A
representation-theoretical interpretation is given in Corollary
\ref{vanderjeugt}. Given $\la\in X_+^{1|n}$, we shall write
\begin{eqnarray} \la^{(0)}=\la, \quad
\la^{(l+1)}=(\la^{(l)})^{\texttt L},\quad l\ge 0.  \label{order}
\end{eqnarray}
The Bruhat order on $X_+^{1|n}$ is the partial order such that for
$\la, \mu \in X_+^{1|n}$, $\mu\prec \la$ if and only if
$\mu=\la^{(l+1)}$ for some $l\ge 0$.

Denote by $Y^{1|n}_+$ the set of the $(n+1)$-tuples $f=(f_{-1}
\mid f_1, f_2, \dots, f_n)$ of integers such that $f_1< f_2<\dots<
f_n<0.$ There is a bijection
\begin{align}
X_+^{1|n}\rightarrow Y^{1|n}_+, \quad \lambda \mapsto f_\lambda,
\label{f-map}
\end{align}
where $f_\la$ is specified by
\begin{align*}
(f_\lambda)_{-1}=(\lambda+\rho, \, \epsilon), \quad
(f_\lambda)_i=(\lambda+\rho, \, \delta_i), \ i\ge 1.
\end{align*}
An element $f_\la \in Y^{1|n}_+$ will also be called atypical
(respectively typical) if $\la \in X_+^{1|n}$ is atypical
(respectively typical). Note that $f \in Y^{1|n}_+$ is atypical if
$|f_{-1}|= -f_i$ for some $i\ge 1$. If $f=f_\lambda$ for some
atypical $\lambda\in X_+^{1|n}$, we set $f^{\texttt L} =
f_{\lambda^{\texttt L}}.$ The Bruhat order on $X_+^{1|n}$ induces
a partial order on $Y^{1|n}_+$ via the bijection (\ref{f-map}).

The description of $f^{\texttt L}$ for a given atypical $f$ is
divided into three cases as follows.

(I) $f_{-1} = f_i <0$ for some $1\le i \le n$. Let $d$ be the
largest integer such that $d<f_i$ and $d\notin \{f_1, \dots,
f_n\}$.  Then, $f^{\texttt L} = (d |f_1, \dots, \hat{f_i}, \dots,
f_n, d)^+$, where $\hat{f_i}$ denotes the removal of $f_i$, and
$+$ denotes the rearrangement of $f_1, \dots, \hat{f_i}, \dots,
f_n, d$ in a decreasing order.

(II) $f_{-1} = - f_i >0$ for some $1\le i \le n$, and $\{f_1, f_2,
\dots, f_n\}$ does not contain $\{-1, -2, \dots, f_i\}$ as a
subset. Let $c$ be the largest integer such that $-1 \ge -c > f_i$
and $-c \notin \{f_1, \dots, f_n\}$. Then, $f^{\texttt L} = (c
|f_1, \dots, \hat{f_i}, \dots, f_n, -c)^+$.

(III) $f_{-1} = - f_i >0$ for some $1\le i \le n$, and the set
$\{f_1, f_2, \dots, f_n\}$ contains $\{-1, -2, \dots, f_i\}$ as a
subset. Then $f^{\texttt L} =(-f_{-1}|f_1, \dots, f_n)$.

\begin{lemma} \label{dualweight}
Let $\beta =2n \epsilon$, the sum of all odd positive roots. Let
$w_0$ be the longest element in $W$. Then, $f_{\beta -w_0 \la}$ is
obtained from $f_\la$ by changing the sign on the (-1)st
component. That is, $f_{\beta -w_0 \la} = -w_0 f_\la.$
\end{lemma}

\begin{proof}
Note that $w_0 \la =\la_{-1} \epsilon -\sum_{i=1}^n \la_i\delta_i$
for $\la =\la_{-1} \epsilon +\sum_{i=1}^n \la_i \delta_i$. The
lemma now follows by unravelling the definitions. \end{proof} It
follows from Lemma~\ref{dualweight} and (I), (II), (III) above
that
\begin{equation} \label{dualdual}
\beta -w_0 \la
 = (\beta -w_0 \la^{\texttt L})^{\texttt L}.
\end{equation}

\subsection{Quantum group}

The quantum group $U_q({\mathfrak g\mathfrak l}(\infty))$ is the
$\mathbb Q(q)$-algebra generated by $E_a, F_a, K^{\pm}_a, a \in
\Z$, subject to the following relations
\begin{eqnarray*}
 K_a K_a^{-1} =K_a^{-1} K_a =1, &&
 K_a K_b = K_b K_a, \\
 K_a E_b K_a^{-1} = q^{\delta_{a,b} -\delta_{a,b+1}} E_b, &&
 K_a F_b K_a^{-1} = q^{\delta_{a,b+1}-\delta_{a,b}}
 F_b, \\
 E_a F_b -F_b E_a &=& \delta_{a,b} (K_{a,a+1}
 -K_{a+1,a})/(q-q^{-1}), \\
 E_a E_b = E_b E_a,   &&
 F_a F_b = F_b F_a,  \qquad\qquad \text{if } |a-b|>1,\\
 E_a^2 E_b +E_b E_a^2 &=& (q+q^{-1}) E_a E_b E_a,  \qquad\qquad \text{if } |a-b|=1, \\
 F_a^2 F_b +F_b F_a^2 &=& (q+q^{-1}) F_a F_b F_a,  \qquad\qquad \text{if } |a-b|=1.
\end{eqnarray*}
Here and further $K_{a,a+1} :=K_aK_{a+1}^{-1}, a\in \Z$. Define
the bar involution on $U_q({\mathfrak g\mathfrak l}(\infty))$ to
be the anti-linear automorphism $-$ such that $\overline{E}_a =
E_a, \overline{F}_a =F_a, \overline{K}_a =K_a^{-1}.$ As usual {\em
anti-linear} means $q \mapsto q^{-1}$. We will sometimes also
write $E_{a}=E_{a,a+1}$ and $F_a=E_{a+1,a}$.

Let $\mathbb V$ be the natural $U_q({\mathfrak g\mathfrak
l}(\infty))$-module with basis $\{v_a\}_{a\in\Z}$ and $\mathbb W
:=\mathbb V^*$ the dual module with basis $\{w_a\}_{a\in\Z}$ such
that
\begin{align}\label{dual}
w_a (v_b) = (-q)^{-a} \delta_{a,b}.
\end{align}
The action of the Chevalley generators on these basis elements are
given explicitly by:
\begin{align*}
K_av_b=q^{\delta_{ab}}v_b,\quad E_av_b=\delta_{a+1,b}v_a,\quad
F_av_b=\delta_{a,b}v_{a+1},\\
K_aw_b=q^{-\delta_{ab}}w_b,\quad E_aw_b=\delta_{a,b}w_{a+1},\quad
F_aw_b=\delta_{a+1,b}w_{a}.
\end{align*}
We shall use the same comultiplication $\Delta$ on $U_q({\mathfrak
g\mathfrak l}(\infty))$ as in \cite{B1,CWZ}:
\begin{eqnarray}
 \Delta (E_a) &=& 1 \otimes E_a + E_a \otimes K_{a+1, a}, \nonumber\\
 \Delta (F_a) &=& F_a \otimes 1 +  K_{a, a+1} \otimes F_a,
\label{comult}\\
 \Delta (K_a) &=& K_a \otimes K_{ a}.\nonumber
\end{eqnarray}

We denote the Iwahori-Hecke algebra of type $A$ by $\mathcal
H_{n}$, which is the $\mathbb Q(q)$-algebra generated by $H_i$,
where $1 \le i \le n-1$, subject to the relations
\begin{eqnarray*}
(H_i -q^{-1})(H_i +q) = 0,\;
 H_i H_{i+1} H_i = H_{i+1} H_i H_{i+1},\;
 H_i H_j = H_j H_i \;(|i-j| >1).
\end{eqnarray*}
For $x \in S_{n}$ with a reduced expression $x=s_{i_1} \cdots
s_{i_r}$, we set
%\begin{equation*}
$H_x :=H_{i_1} \cdots H_{i_r}.$
%\end{equation*}
The bar involution $-$ on $\mathcal H_{n}$ is the unique
anti-linear automorphism defined by $\overline{H_x}
=H_{x^{-1}}^{-1}$ for all $x \in S_{n}$. We let
$H_0 := \sum_{x \in S_{n}} (-q)^{\ell(x) - \ell (\sigma_0)} H_x$
where $\sigma_0$ is the longest element in $S_n$.

\section{Canonical and dual canonical bases on a Fock space}
\label{Fock-space}
\subsection{The Fock space}

Denote by $\mathbb V_+$, $\mathbb V_0$, and $\mathbb V_-$ the
subspaces of $\mathbb V$ spanned by $v_i$ with $i>0$, $i=0$, and
$i<0$, respectively. Denote by $\mathcal U$ (respectively
$\mathcal U_+$) the subalgebras of $U_q({\mathfrak g\mathfrak
l}(\infty))$ generated by $E_i, F_i, K_{i+1},K^{-1}_{i+1}$, $i \le
-2$ (respectively by $E_i , F_i, K_{i}, K^{-1}_i$, $i \ge 1$).
Then $\mathbb V_-$ (respectively $\mathbb V_+$) is the natural
$\mathcal U$ (respectively $\mathcal U_+$) module.

\begin{lemma}
\begin{enumerate}
\item There is an algebra isomorphism $\phi: \mathcal U_+
\stackrel{\sim}{\mapsto} \mathcal U$ given by
$$E_{i} \stackrel{\phi}{\mapsto} F_{-i-1}, \quad
F_{i} \stackrel{\phi}{\mapsto} E_{-i-1}, \quad K^{\pm 1}_i
\stackrel{\phi}{\mapsto} K^{\pm 1}_{-i} \quad (i>0).$$
%
%or in a more symmetric form,
%%
%$$E_{i, i+1} \stackrel{\phi}{\mapsto} E_{-i,-i-1}, \quad
%E_{i+1, i} \stackrel{\phi}{\mapsto} E_{-i-1,-i}, \quad K_{i}
%\stackrel{\phi}{\mapsto} K_{-i} \quad (i>0).$$

\item
The composition of $\phi^{-1}$ with the natural action of
$\mathcal U_+$ on $\mathbb V_+$ defines a $\mathcal U$-module
structure on $\mathbb V_+$.  Furthermore the linear map $\mathbb
V_+ \rightarrow \mathbb V_-$, $v_i \mapsto v_{-i}$, $i>0$, is an
isomorphism of $\mathcal U$-modules.
\end{enumerate}
\end{lemma}
\begin{proof}
Part (1) follows by checking the defining relations, while (2)
follows from the explicit formulae of the actions of $\mathcal U$
on $\mathbb V_-$ and $\mathcal U_+$ on $\mathbb V_+$. \end{proof}

Denote by $\mathbb W_-$ the $\mathcal U$-module which is dual to
the natural $\mathcal U$-module $\mathbb V_-$, with generators
$w_i$, $i<0$, normalized as in (\ref{dual}). Define the
$\cU$-module $\bigotimes^n\mathbb W_-$  via the usual
comultiplication
$\Delta^{(n-1)}=(id^{\otimes(n-2)}\otimes\Delta)\dots
(id\otimes\Delta)\Delta$. By the Schur-Jimbo duality, $\mathcal
H_n$ acts on $\bigotimes^n\mathbb W_-$ and this action commutes
with the action of $\mathcal U$.  Define $\bigwedge^n \mathbb W_-$
to be the quotient of $\bigotimes^n\mathbb W_-$ by the kernel of
$H_0$. Denote the image of $w_{f_1}\otimes \cdots \otimes w_{f_n}$
in $\bigwedge^n \mathbb W_-$ by $w_{f_1}\wedge \cdots \wedge
w_{f_n}$, for $w_{f_1}, \ldots, w_{f_n} \in \mathbb W_-$.

Consider the following Fock space
$$\mathcal F := \mathbb V \otimes \bigwedge^n \mathbb
W_-.$$
For $f =(f_{-1} | f_1, \cdots, f_n)\in Y^{1|n}_+$, let
\begin{eqnarray}
K_f &:=&v_{f_{-1}}\otimes w_{f_1}\wedge w_{f_2}\wedge \cdots
\wedge w_{f_n}.
\end{eqnarray}
Then the $K_f$ form a basis of $\cF$.

Define an action of the quantum group $\mathcal U$ on the space
$\cF$ in the following way. Let $\mathcal F_\bullet :=\mathbb
V_\bullet \otimes \bigwedge^n \mathbb W_-$ for $\bullet = +, -,
0.$ The action of $\mathcal U$ on $\mathcal F_-$ is defined
exactly as in (\ref{comult}) via $\Delta$.
%\begin{eqnarray*}
%\Delta(E_{-i-1, -i}) &=& E_{-i-1, -i}\otimes K_{-i,-i-1}
%+ 1\otimes E_{-i-1, -i}, \\
%\Delta(E_{-i, -i-1}) &=& E_{-i, -i-1}\otimes 1 + K_{-i-1,-i}
%\otimesE_{-i, -i-1}, \\
%\Delta(K_{-i})&=&K_{-i}\otimes K_{-i}.
%\end{eqnarray*}
The action of $\mathcal U$ on $\mathcal F_+$ is via $\Delta'
=(1\otimes \phi) \circ \Delta \circ\phi^{-1}$, which is a mixture
of the comultiplication $\Delta$ and the isomorphism $\phi$: for
$i>0$,
\begin{eqnarray*}
\Delta'(E_{-i,-i-1}) &=& 1\otimes E_{-i,-i-1} + E_{i,i+1} \otimes
K_{-i-1,-i}, \\
\Delta'(E_{-i-1,-i}) &=& E_{i+1,i} \otimes 1 + K_{i,i+1} \otimes
E_{-i-1,-i}, \\
\Delta'(K_{-i}) &=& K_{i} \otimes K_{-i}.
\end{eqnarray*}
The action of $\cU$ on $\mathcal F_0$ is defined by $x \mapsto 1
\otimes x$ for every $x \in \cU$
%$E_{-i-1, -i} \mapsto 1\otimes E_{-i-1, -i}, E_{-i, -i-1} \mapsto
%$1 \otimes E_{-i, -i-1},$ and $K_{-i}\mapsto 1 \otimes K_{-i}$,
which is compatible with either $\Delta$ or $\Delta'$. Putting
together, we have defined an action of $\cU$ on $\mathcal F
=\mathcal F_+ \oplus \mathcal F_0  \oplus \mathcal F_-.$

We also define the following operators on $\mathcal F$:
$$E_{-1} := E_{-1,0} \otimes (K_{-1})^{-1} + E_{1,0} \otimes 1,
\qquad F_{-1} := E_{0,-1} \otimes 1 + E_{0,1} \otimes K_{-1}.$$
%
%Formally $e$ is $\Delta(E_{-1,0}) +\Delta' (E_{-1,0}),$ and
%similarly for $f$.
The following lemmas can be proved by straightforward
calculations.

\begin{lemma} \label{ef}
For $g=(g_{-1}|g_1,  \dots, g_n)\in Y_+^{1|n}$, let
$g^\pm=(g_{-1}\pm 1|g_1, \dots, g_n).$   The actions of $E_{-1}$
and $F_{-1}$ on $K_g$ vanish unless $g_{-1}=0, \pm 1$, or
$g_n=-1$. In these cases,
\begin{enumerate}
\item if $g_{n}=-1,$ then
\begin{eqnarray*}
E_{-1} (K_g)= K_{g^+} +q K_{g^-},& \quad \text{if} \ g_{-1}=0, \\
F_{-1} (K_g) =  q^{-1} K_{g^-}, & \quad \text{if} \ g_{-1}=1,\\
F_{-1} (K_g) =  K_{g^+}, & \quad \text{if} \ g_{-1}=-1;
\end{eqnarray*}
\item if $g_{n}\ne -1,$ then
\begin{eqnarray*}
E_{-1} (K_g)=K_{g^+} + K_{g^-}, & \quad \text{if} \ g_{-1}=0,\\
F_{-1} (K_{g}) =   K_{g^-}, & \quad \text{if} \ g_{-1}=1,\\
F_{-1} (K_{g}) = K_{g^+}, & \quad \text{if} \ g_{-1}=-1.
\end{eqnarray*}
\end{enumerate}
\end{lemma}

%\begin{lemma} \label{ef-fe}
%Keep the notation of Lemma \ref{ef}. Set $h =ef -fe$. Then,
%\begin{enumerate}
%\item for $g=(0|g_1, \dots, g_{n-1}, -1)$, we have
% $h (K_g)
% = -(q +q^{-1}) K_g.
% $
%\item for $g=(0|g_1, \dots, g_{n})$ with $g_n\ne -1$, we have
%\begin{eqnarray*}
% h (K_g) = -2 K_g,&&
% h (K_{g^+} +q K_{g^-})
% = 2 (K_{g^+} +q K_{g^-}).
%\end{eqnarray*}
%\end{enumerate}
% \end{lemma}
{From} now on, by abuse of notation, and we will refer to $E_a,
F_a (a \ge -1)$ as the {\em Chevalley generators}.

\subsection{The canonical and dual canonical bases}

\begin{proposition} \label{procedure}
For every atypical $f\in Y_+^{1|n}$, there exist a typical $g\in
Y_+^{1|n}$ and a sequence of Chevalley generators $X_1, \dots,
X_r$ such that
\[
X_1\cdots X_r (K_g) =K_f +q K_{f^{\texttt L}}.\]
\end{proposition}
\begin{proof}
We explicitly construct the sequence of Chevalley generators $X_1,
\dots, X_r$ and a typical weight $g$ for every atypical $f$. There
are three cases to consider according to Subsection~\ref{bruhat}.

(I) $f_{-1}=f_i$ for a fixed $i \ge 1$. Set $k=-f_{-1}$. There
exists an $l\ge 0$ such that
\begin{eqnarray}  f &=&(-k | f_1,\cdots,
f_{i-l-1}, -k-l, -k-l+1, \dots, -k, f_{i+1}, \cdots, f_n), \label{case1} \\
f^{\texttt L}&=& (-k-l-1 | f_1, \cdots, f_{i-l-1}, -k-l-1, -k-l,
\dots, -k-1, f_{i+1}, \cdots, f_n),\nonumber
\end{eqnarray}
with $f_{i-l-1}<  -k-l-1$. Then,
\[
%E_{-k-l-1, -k-l} E_{-k-l, -k-l+1} \cdots E_{ -k-1,-k}
E_{-k-l-1} E_{-k-l} \cdots E_{ -k-1} \left(K_g\right) =K_f +q
K_{f^{\texttt L}}.
\]
with the typical element
\[ g := (-k| f_1,\cdots, f_{i-l-1}, -k-l-1, -k-l, \dots,
-k-1, f_{i+1}, \cdots, f_n). \]

(II) $f_{-1}=-f_i$ for some $i \ge 1$, and $f^{\texttt L}_{-1}>0$.
There is an $l\ge 0$ such that
\begin{eqnarray}
f &=&(k+l+1 | f_1, \cdots, f_{i-1}, -k-l-1, -k-l, \dots, -k-1,
f_{i+l+1},
\cdots, f_n), \nonumber \\
f^{\texttt L}&=& (k| f_1, \cdots, f_{i-1}, -k-l, -k-l+1, \dots,
-k, f_{i+l+1}, \cdots, f_n), \label{case2} \nonumber
\end{eqnarray}
with $f_{i+l+1}>-k$, where $k=f^{\texttt L}_{-1}$. Then,
\[
%E_{-k, -k-1} E_{-k-1, -k-2} \cdots E_{-k-l,-k-l-1}
F_{-k-1} F_{-k-2} \cdots F_{-k-l-1} \left(K_g\right) =K_f +q
K_{f^{\texttt L}}
\]
with the typical element
\[ g := (k+l+1| f_1, \cdots, f_{i-1}, -k-l, -k-l+1, \dots, -k,
f_{i+l+1}, \cdots, f_n).\]

(III) $f_{-1}=-f_i$ for some $i \ge 1$, and $f^{\texttt
L}_{-1}<0$. Then
\begin{eqnarray}
f &=&(k| f_1, \dots, f_{n-k}, -k, -k+1, \dots, -1),\nonumber \\
f^{\texttt L}&=& (-k| f_1, \dots, f_{n-k}, -k, -k+1, \dots, -1),
\label{case3} \nonumber
\end{eqnarray}
where $k=f_{-1}$. Then
$$
%E_{-k, -k+1}  \dots   E_{-2, -1} E_{-1, 0}
E_{-k}  \dots   E_{-2} E_{-1} (K_g) =K_f +q K_{f^{\texttt L}}$$
with the typical element
$g :=(0| f_1, \dots, f_{n-k}, -k, -k+1, \dots, -1).$ \end{proof}

%The operators $E_{-a}$, $F_{-a}$, $a\ge 1$, will be referred to as
%the Chevalley generators.
We define a bar-involution $-$ on $\cU$ by declaring that it fixes
all the Chevalley generators and sends $K_{-i}$ to $K_{-i}^{-1}$.
\begin{theorem}
\begin{enumerate}
\item
There exists a unique anti-linear bar involution $-$ on a suitable
completion $\widehat{\mathcal F}$ of ${\mathcal F}$ such that
\begin{enumerate}
\item[(i)] $\overline{K_f} = K_f$ for all typical $f \in
Y_+^{1|n}$;

\item[(ii)]
$\overline{X u} = \overline{X} \overline{u}$, $\overline{E_{-1} u}
= E_{-1} \overline{u}$ and $\overline{F_{-1} u} = F_{-1}
\overline{u}$, for all $X \in \mathcal U$ and $u \in
\widehat{\mathcal F}$.
\end{enumerate}
\item
There exists unique {\em canonical basis} $\{U_f\}$ and {\em dual
canonical basis} $ \{L_f \}$, where ${f \in Y_+^{1|n}}$,  for
$\widehat{\mathcal F}$ such that
\begin{enumerate}
\item[(i)] $\overline{U}_f =U_f$ and $\overline{L}_f =L_f$;
\item[(ii)] $U_f \in K_f + \widehat{\sum}_{g \prec f} q \Z [q] K_g$
 and $L_f \in K_f + \widehat{\sum}_{g \prec f} q^{-1} \Z [q^{-1}]
 K_g$.
\end{enumerate}
\item
$U_f =L_f =K_f$ for typical $f\in Y_+^{1|n}$. For every atypical
$f\in Y_+^{1|n}$, we have
\begin{eqnarray}
U_f =K_f +q K_{f^{\texttt L}}, \qquad
L_f = K_f +\sum_{l=1}^\infty (-q^{-1})^l K_{f^{(l)}}
\label{irreps}
\end{eqnarray}
where $f^{(1)}=f^{\texttt L}$ and $f^{(l+1)}=(f^{(l)})^{\texttt
L}$.
\end{enumerate}
\end{theorem}

\begin{proof}
Proposition~\ref{procedure} and the requirement (ii) of the bar
map imply that $K_f +q K_{f^{\texttt L}}$ for every atypical $f$
is bar-invariant. Thus, $K_f +q K_{f^{\texttt L}}$ for all
atypical $f\in Y_+^{1|n}$ together with $K_f$  for all typical
$f\in Y_+^{1|n}$ form a bar-invariant basis of $\widehat{\mathcal
F}$. This proves the uniqueness of the bar map.

Since $\overline{K_f +q K_{f^{\texttt L}}} =K_f +q K_{f^{\texttt
L}}$ for $f$ atypical, we obtain
\begin{eqnarray}
\overline{K_f} =K_f +q K_{f^{\texttt L}} -q^{-1}
\overline{K_{f^{\texttt L}}}. \label{K-bar}
\end{eqnarray}
By iterating the relation (\ref{K-bar}) we obtain that
\begin{eqnarray*}
\overline{K_f} = K_f+(q-q^{-1})\sum_{i=1}^\infty
(-q)^{1-i}K_{f^{(i)}}.
\end{eqnarray*}
It follows that the bar map is indeed an involution with the
property that $\overline{K_f}$ equals $K_f$ plus lower terms in
Bruhat order for every $f\in Y_+^{1|n}$. The existence and
uniqueness of the canonical and dual canonical bases now follows
routinely from the bar involution with such a property \cite{KL}.

Clearly, for every typical $f\in Y_+^{1|n}$, we have $U_f =L_f
=K_f$. By the uniqueness of the canonical basis, $U_f =K_f +q
K_{f^{\texttt L}}$ for $f$ atypical. Denote the RHS of
(\ref{irreps}) by $\mathfrak L_f$. It follows from (\ref{K-bar})
that $ \overline{\mathfrak L_f} = K_f -q^{-1} \overline{\mathfrak
L_{f^{\texttt L}}}.$ Iterating this relation we obtain
$\overline{\mathfrak L_f} = \mathfrak L_f$. Thus $L_f =\mathfrak
L_f$ by the uniqueness of the dual canonical basis.

It remains to check that the bar map on $\mathcal F$ indeed
satisfies the compatibility condition (ii) of (1). For the
generators $E_{-1}$ and $F_{-1}$, this follows from
Lemma~\ref{ef}.

Now consider the Chevalley generators of $\cU$. This requires a
tedious (albeit elementary) case by case verification that
$\overline{X_{-a}(U_f)}=X_{-a}(U_f)$, for all $f\in Y_+^{1|n}$ and
all Chevalley generators $X_{-a}$ of $\cU$. If $f\in Y_+^{1|n}$ is
typical, then for any Chevalley generator $X_{-a}$ of $\cU$,
$X_{-a}(U_f)$ is either zero or equal to some $U_g$. This can be
established by separately analyzing the two cases with
$|f_1|+f_i\ne \pm 1$ for all $i>0$  and  with $|f_1|+f_i= \pm 1$
for some $i>0$ respectively.

We divide the atypical elements of $Y_+^{1|n}$ into three cases as
in the proof of Proposition \ref{procedure}. (I) $f$ and
$f^{\texttt{L}}$ are given by equation (\ref{case1}). If $l\ge 1$,
we can show that for all the Chevalley generators $X_{-a}$ of
$\cU$, $X_{-a}(U_f)$ is either zero or equal to some $U_g$. This
is also true for all the $X_{-a}$ but $E_{-k-1}$ when $l=0$. In
the latter case, $E_{-k-1}(U_f) = (q+q^{-1}) U_g$ with $g= (-k-1|
f_1, \dots, f_{i-1}, -k, f_{i+1}, \dots, f_n)$. The case (II) is
analogous, thus we omit the details. In the case (III),
$X_{-a}(U_f)$ is either zero or equal to some $U_g$ if $X_{a}\ne
F_{-k}$. When $X_{a}= F_{-k}$, we have
\[ F_{-k}(U_f) = \left\{\begin{array}{l l} 0, & \text{if} \  f_{n-k}=
-k-1,\\
U_{g_+} + U_{g_-}, & \text{if} \  f_{n-k}\ne -k-1,
\end{array}\right.\]
where $g_\pm = (\pm k| f_1, \dots, f_{n-k}, -k-1, -k+1, \dots,
-1)$. This completes the proof of the theorem. \end{proof}

\section{Representation theory of ${\mathfrak g}$}
\subsection{Characters of the tilting and irreducible ${\mathfrak g}$-modules}

Let $\sum_{i=0}^\infty \Z\epsilon_i$ denote the free abelian group
with basis $\epsilon_i$, $i=0, -1, -2, \dots$. We define a map
$\text{wt}: Y^{1|n}_+\mapsto \sum_{i=0}^\infty\Z\epsilon_i$ by
\begin{eqnarray*}
f=(f_{-1} \mid f_1, f_2, \dots, f_n)\mapsto
\text{wt}(f)=\epsilon_{-|f_{-1}|} - \epsilon_{f_{1}} -
\epsilon_{f_{2}}- \cdots - \epsilon_{f_{n}}.
\end{eqnarray*}
It was stated in \cite{K2} and proved in \cite[Theorem~1.2]{Pe}
that $\la, \mu\in X_+^{1|n}$ correspond to the same central
character only if $\text{wt}(f_\la)=\text{wt}(f_\mu)$. On the
other hand, exactly when $\la$ and $\mu$ are atypical,
$\text{wt}(f_\la)=\text{wt}(f_\mu)$ implies that $\la$ and $\mu$
are in the same block. The block corresponding to $\la$ will be
denoted by $\cO^+_\la$. Evidently $\cO^+$ is a direct sum of
blocks corresponding to different central characters.

Given a $\mathfrak g$-module $M$ in the category $\cO^+$, we endow
the dual $M^*$ with the usual $\mathfrak g$-module structure.
Further twisting the $\mathfrak g$-action on $M^*$ with the
automorphism of $\mathfrak g$ given by $e_i \mapsto -f_i, f_i
\mapsto -e_i, h_i \mapsto h_i$ for $1\le i\le n$ and $e_0 \mapsto
f_0, f_0 \mapsto -e_0, h_0 \mapsto h_0$, we obtain another
$\mathfrak g$-module denoted by $M^\tau$. We have $(M^\tau)^\tau
\cong M$.

We shall consider translation functors on the category $\cO^+$.
For any $M\in \cO^+$ belonging to the block of a weight $\la \in
X_+^{1|n}$, we define for $a=1, 2, \dots$
\begin{eqnarray*}
\tilde{E}_{-a}(M) &=& \text{pr}_{\text{wt}
(f_\la)+\epsilon_{-a}-\epsilon_{-a+1}}(\C^{2|2n}\otimes
M),\\
\tilde{F}_{-a}(M) &=& \text{pr}_{\text{wt}
(f_\la)+\epsilon_{-a+1}-\epsilon_{-a}}(\C^{2|2n}\otimes M).
\end{eqnarray*}
Here for $\gamma\in\sum_{i=0}^\infty\Z\epsilon_i$,
$\text{pr}_{\gamma}:\cO^+\rightarrow
\bigoplus_{\text{wt}(\mu)=\gamma}\cO^+_\mu$ stands for the
canonical projection. Such functors are exact and their left and
right adjoints are of the same form.

Let $K(\cO^+)$ denote the Grothendieck group of $\cO^+$. For
$M\in\cO^+$ the expression $[M]$ denotes the corresponding element
in $K(\cO^+)$. We shall use the same notation to denote the
operators on $K(\cO^+)$ corresponding to $\tilde{E}_{-a}$ and
$\tilde{F}_{-a}$ respectively. By checking the tensor product of
$K(\la)$ with the natural module $\C^{2|2n}$, we can easily prove
the following result.

\begin{proposition}  \label{translation}
The linear map $j: K(\mathcal O^+) \rightarrow \mathcal F|_{q=1}$,
$K(\la) \mapsto K_{f_\la}$, is an isomorphism of vector spaces.
Furthermore, for $a \le -1$ we have
$$E_{-a} j(-) =j(\tilde{E}_{-a}(-)), \quad F_{-a} j(-)
=j(\tilde{F}_{-a}(-)).$$
\end{proposition}

We say that an object $M\in\cO^+$ has a {\em Kac flag}, if $M$ has
a filtration of submodules from $\cO^+$ such that each successive
quotient is isomorphic to some Kac module. The general theory of
finite dimensional tilting modules as explained in \cite{So, B2}
applies to the Lie superalgebra ${\mathfrak g}$ as well.  We
denote by $U(\la)\in \cO^+$ the {\em tilting module} associated
with $\la\in X_+^{1|n}$.  It is the unique indecomposable object
in $\cO^+$ satisfying: (1) $U(\la)$ has a Kac flag with $K(\la)$
at the bottom; (2) $\text{Ext}^1_{\cO^+}(K(\mu),U(\la))=0$ for all
$\mu\in X_+^{1|n}$. Denote by $(U(\la): K(\mu))$ the multiplicity
of the Kac module $K(\mu)$ in a Kac flag of $U(\la)$, and by
$[K(\mu): L(\nu)]$ the multiplicity of $L(\nu)$ in a composition
series of $K(\mu)$. Recall $\beta =2n \epsilon$. Following
\cite{B1, B2}, we have
\begin{eqnarray}
K(\la)^* \cong K(\beta -w_0 \la); \quad U(\la)^* \cong U(\beta
-w_0 \la). \nonumber\\
(U(\la): K(\mu)) = [K(\beta -w_0 \mu): L(\beta -w_0\la)].
\label{betadual}
\end{eqnarray}

Below we have the following analogue of Theorem~4.37 in \cite{B1}.
\begin{theorem} \label{main}
Let $\la\in X_{1|n}^+$. Then,
\begin{enumerate}
\item If $\la$ is typical, then $U(\la) =K(\la) =L(\la)$.
If $\la$ is atypical, there exist a sequence of translation
functors $X_1, \dots, X_r$ and a typical $\mu$ in $X_{1|n}^+$ such
that $U(\la)=X_1 \cdots X_r U(\mu)$. Furthermore, $U(\la)$  has
the following $2$-step Kac flag: $0 \rightarrow K(\la) \rightarrow
U(\la) \rightarrow K(\la^{\texttt L}) \rightarrow 0$.

\item $j ([U(\la)]) = U_{f_\la}|_{q=1}.$

\item $U(\la)$ is the projective cover of $L(\la^{\texttt L})$.

\item $U(\la) \cong U(\la)^\tau.$
\end{enumerate}
\end{theorem}

\begin{proof}
Let us first assume the validity of (1). Part (2) immediately
follows from Proposition \ref{translation} and Proposition
\ref{procedure}. Part (1) implies that there is an epimorphism
$U(\la)\rightarrow L(\la^L)$ and $U(\la)$ is indecomposable.  We
have $\text{Hom}_{\cO^+}(U(\la),M)\cong
\text{Hom}_{\cO^+}(K(\mu),Y_r \cdots Y_1 M)$, where $Y_a$ is the
translation functor corresponding to the adjoint functor of $X_a$.
Thus $U(\la)$ is projective and (3) holds. Part (4) also follows
readily from (1) by an induction argument, since $\tau$ commutes
with the translation functors.

So it remains to prove (1). The typical case is clear. Now let us
fix an atypical $\la$. Then by Proposition \ref{procedure} there
exists a typical $\nu$ and a sequence of Chevalley generators
$X_1, \dots, X_r$ such that $U_{f_\la} =X_1\cdots X_r
(U_{f_\nu})$. By abuse of notation, we shall also denote by $X_1,
\dots, X_r$ the corresponding translation functors.

Clearly (1) holds for a typical $\nu$. By
Proposition~\ref{translation}, we have
$$j[X_1\cdots X_r  U(\nu)] =X_1\cdots X_r  U_{f_\nu}(1)
=U_{f_\la}(1).$$
This and the formula for the canonical basis element $U_{f_\la}
=K_{f_\la} +qK_{f_{\la^{\texttt L}}}$ imply the following identity
in $K(\cO^+)$:
\begin{equation}  \label{twoterm}
[X_1\cdots X_r  U(\nu)] =[K(\la)] +[K (\la^{\texttt L})].
\end{equation}
Now it is easy to see that if $M\in\cO^+$ has a Kac flag, then the
translation functor applied to $M$ produces a module with a Kac
flag, which in turn implies that any direct summand of it also has
a Kac flag. Since $\text{Ext}^1_{\cO^+}(K(\mu),X_1\cdots X_r
U(\nu))=\text{Ext}^1_{\cO^+}(Y_r\cdots Y_1K(\mu), U(\nu))=0$, we
have $\text{Ext}^1_{\cO^+}(K(\mu),U)=0$ for any direct summand $U$
of $X_1\cdots X_r U(\nu)$.  Hence by (\ref{twoterm}) $X_1\cdots
X_r  U(\nu)$ is a direct sum of tilting modules and contains
$U(\la)$ as a direct summand.

If we can show that $K (\la^{\texttt L})$ (besides the obvious one
$K(\la)$) appears in a Kac flag for $U(\la)$, then by
(\ref{twoterm}) again there will be no more tilting module as a
direct summand of $X_1\cdots X_r  U(\nu)$ and we will be done.
This latter claim is indeed true, since by (\ref{dualdual}) and
(\ref{betadual}),
\begin{eqnarray*}
(U(\la): K (\la^{\texttt L}))
 &=& [K (\beta -w_0 \la^{\texttt L}): L(\beta -w_0 \la)] \\
 &=& [K (\beta -w_0 \la^{\texttt L}):
 L((\beta -w_0 \la^{\texttt L})^{\texttt L})] \ge 1.
 \end{eqnarray*}
The last inequality $[K(\mu): L(\mu^{\texttt L})] \ge 1$ for every
atypical $\mu$ will be established in Lemma \ref{submodule} below,
independently of van der Jeugt's theorem \cite{V}. \end{proof}

Using (\ref{dualdual}), (\ref{betadual}) and Theorem~\ref{main},
we have obtained a new proof of van der Jeugt's main theorem.
\begin{corollary}\label{vanderjeugt} \cite{V}
 For $\la\in X^{1|n}_+$ atypical, there is a short exact sequence of
 ${\mathfrak g}$-modules
 $$ 0 \rightarrow L(\la^{\texttt L}) \rightarrow K(\la)
 \rightarrow L(\la) \rightarrow 0. $$
\end{corollary}
Let $\rho_0 =\frac12 \sum_{\alpha \in \Delta_0^+}\alpha$. Write
$$D_0:=\sum_{w\in W}(-1)^{l(w)}w(e^{\rho_0}).$$
\begin{corollary}\label{BL} \cite{V}
let $\la \in X_+^{1|n}$ be atypical weight with $(\la +\rho,
\gamma)=0$ for some $\gamma \in \Delta_1^+$. Then
\begin{eqnarray} \label{BLformula}
\ch L(\la) =\frac1{D_0} \sum_{w\in W} (-1)^{l(w)} w\left(
e^{\la+\rho_0} \prod_{\alpha \in \Delta_1^+\backslash\{\gamma\}}
(1+e^{-\alpha})\right).
\end{eqnarray}
\end{corollary}
\begin{proof}
This follows from Corollary~\ref{vanderjeugt} and the fact that
\begin{align*}
\text{RHS of } (\ref{BLformula}) &=\frac{\prod_{\alpha \in
\Delta_1^+} (1+e^{-\alpha})}{D_0} \sum_{w\in W} (-1)^{l(w)}
w\left(
e^{\la+\rho_0} /(1+e^{-\gamma})\right) \\
&=\frac{\prod_{\alpha \in \Delta_1^+}
(1+e^{-\alpha})}{\prod_{\alpha\in\Delta_0^+} (1-e^\alpha)}
\sum_{k\ge 0}(-1)^k\sum_{w\in W} (-1)^{l(w)} w\left(
e^{\la+\rho_0-k\gamma} \right)e^{-\rho_0} \\
&=\sum_{k\ge 0}(-1)^k\prod_{\alpha \in \Delta_1^+}
(1+e^{-\alpha})\frac{\sum_{w\in W} (-1)^{l(w)}
e^{w(\la+\rho_0-k\gamma)-\rho_0} }{\prod_{\alpha\in\Delta_0^+}
(1-e^\alpha)}
\\
&=\sum_{i \ge 0} (-1)^i \ch K(\la^{(i)}).
\end{align*}
In the last identity we have used the Weyl character formula and
the fact that $\rho_0-\rho$ is $W$-invariant. \end{proof}

\begin{remark}
In \cite{Zou}, $\Ext^i_{\cO^+}(K(\mu), L(\la))$ was computed
explicitly. Denote the dual canonical basis element  $L_f = \sum_g
l_{g f}(q) K_g$. These Kazhdan-Lusztig polynomials $l_{gf}(q)$
have been computed in (\ref{irreps}) for atypical $f$, and $l_{g
f}(q) =\delta_{g,f}$ for typical $f$. Comparing with \cite{Zou},
we see the Serganova-Zou's Kazhdan-Lusztig polynomials coincide
with ours:
$l_{f_\mu f_\la}(-q^{-1})= \sum_{i\ge 0} q^i
\Ext^i_{\cO^+}(K(\mu), L(\la)).$
{From} the theory of highest weight categories of Cline, Parshall
and Scott (cf. \cite[4-f]{B1} for adaptation to superalgebras), we
have
$$\sum_{i \ge 0} \dim \Ext^i_{\cO^+}(L(\mu), L(\la)) q^i=
\sum_{\nu \in X^+(1|n)} \ell_{\nu\mu} (-q^{-1}) \ell_{\nu\la}
(-q^{-1}).$$
\end{remark}
\subsection{A technical lemma}

The following was used in the proof of Theorem \ref{main}.
\begin{lemma}\label{submodule}
$[K(\la): L(\la^{\texttt L})] \ge 1$ for every atypical $\la$.
\end{lemma}
\begin{proof}
Assume that $\la=(\la_{-1}|\la_1,\cdots,\la_n)$ is atypical. There
are two possibilities: $(\la_{-1}-n) +(\la_i+n-i+1)=0$, or
$\la_{-1}-n=\la_i+n-i+1$ for some $i$. We will treat in detail
below the first case when
\begin{equation} \label{zero}
\la_{-1} +\la_i -i+1=0
\end{equation}
and leave the other similar case to the reader.

Let $T_-$ be the product of all odd negative root vectors and let
$v_\la$ be a highest weight vector of the Kac module $K(\la)$.
Then the vector $T_-v_\la$ has weight
$(\la_{-1}-2n|\la_1,\cdots,\la_n).$ Note that $T_-v_\la$ is
highest weight with respect to the Borel subalgebra containing the
same even part but the opposite odd part of the standard Borel
$\mathfrak b$. We apply now odd reflections in the following order
to get back to the standard Borel:
\begin{equation*}
\epsilon+\delta_1,\epsilon+\delta_2,\cdots,\epsilon+\delta_n,
\epsilon-\delta_n,\epsilon-\delta_{n-1},\cdots,\epsilon-\delta_1.
\end{equation*}
Here the usual rule of odd reflection is that if $(\mu,\alpha)=0$
then the highest weight vector is unchanged, and if
$(\mu,\alpha)\not=0$ then the highest weight vector is obtained by
applying the positive root vector corresponding to $\alpha$ to the
previous highest weight vector (cf. for example, \cite{PS}).

Note that (\ref{zero}) implies that $\la_{-1} -n<0$ and thus
$\la_{-1}-2n \not=\la_1$. So after the first step the weight is
$(\la_{-1}-2n+1|\la_1+1,\cdots,\la_n).$ Repeating the process with
the first $n$ odd roots, we end up with the weight
$(\la_{-1}-n|\la_1+1,\cdots,\la_n+1).$ We continue by using now
the odd root $\epsilon-\delta_n$. If $\la_{-1}-n+\la_n+1=0$, then
$i=n$. So if $n\not=i$, then we need to add $\epsilon-\delta_n$
and get $(\la_{-1}-n+1|\la_1+1,\cdots,\la_{n-1}+1,\la_n).$
%Repeating the whole thing now with $\epsilon-\delta_{n-1}$ etc. we
%finally reached
%\begin{equation*}
%(\la_{-1}-i|\la_1+1,\cdots,\la_i+1,\la_{i+1},\cdots,\la_n).
%\end{equation*}
%Evidently $\la_{-1}-i+\la_i +1=0$, so we do not add
%$\epsilon-\delta_i$ to the weight above.  Similarly we do not add
%$\epsilon-\delta_{i-1},\cdots,\epsilon-\delta_j$.
%
%After that we need to add
%$\epsilon-\delta_{1},\cdots,\epsilon-\delta_{j-1}$.
Finally we end up with the weight
\begin{equation*}
\la^{\texttt L} =(\la_{-1}-i+j-1|\la_1,\cdots,\la_{j-1},\la_j+1,
\cdots,\la_i+1,\la_{i+1},\cdots,\la_n).
\end{equation*}
Here $j$ is determined by that $\la_j=\la_{j+1}=\cdots=\la_i$ and
$\la_{j-1}>\la_i$. Note that in the process we did not add
$\epsilon-\delta_i, \epsilon-\delta_{i-1}, \ldots,
\epsilon-\delta_j$ since $\la_{-1}-i+\la_i +1=0$. In this way, we
have obtained a highest weight vector (relative to $\mathfrak b$)
of highest weight $\la^{\texttt L}$.

\end{proof} %%
\section{The composition factors of symmetric tensors}
\label{sec:tensor}

Let $x,\bar{x}$ be $2$ even variables and
$\xi_1,\cdots,\xi_n,\bar{\xi}_1,\cdots,\bar{\xi}_n$ be $2n$ odd
variables. If we let $\C^{2|2n}$ stand for the standard
representation of $\osp(2|2n)$, then we may identify the symmetric
algebra $S(\C^{2|2n})$ with $\C[x,\bar{x},\xi_i,\bar{\xi}_i]$, the
polynomial algebra in the variables $x,\bar{x}$ and
$\xi_1,\cdots,\xi_n,\bar{\xi}_1,\cdots,\bar{\xi}_n$. In this
identification the action of ${\mathfrak g}$ gets identified with
the action of certain linear differential operators whose explicit
formulas are easily written down. The positive simple root vectors
$e_0,e_1,\cdots,e_n$ and the negative simple root vectors
$f_0,f_1,\cdots,f_n$ are:
\begin{align*}
&e_0= x\frac{\partial}{\partial\xi_1}
+\bar{\xi}_1\frac{\partial}{\partial\bar{x}},\quad
f_0 = \xi_1\frac{\partial}{\partial x} -\bar{x}\frac{\partial}{\partial\bar{\xi}_1},\\
&e_i= \xi_i\frac{\partial}{\partial\xi_{i+1}}
-\bar{\xi}_{i+1}\frac{\partial} {\partial\bar{\xi}_i},\quad f_i=
\xi_{i+1}\frac{\partial}{\partial\xi_{i}}
-\bar{\xi}_{i}\frac{\partial} {\partial\bar{\xi}_{i+1}},
\quad i=1,\cdots,n-1,\\
&e_n= \xi_n\frac{\partial}{\partial\bar{\xi}_n},\quad f_n=
\bar{\xi}_n\frac{\partial}{\partial{\xi}_n}.
\end{align*}
By declaring all the variables to have degree $1$ the algebra
$\C[x,\bar{x},\xi_i,\bar{\xi}_i]$ acquires a $\Z$-grading
\begin{align*}
\C[x,\bar{x},\xi_i,\bar{\xi}_i]=\bigoplus_{j=0}^\infty\C[x,\bar{x},\xi_i,\bar{\xi}_i]^j
\cong \bigoplus_{j=0}^\infty S^j(\C^{2|2n}).
\end{align*}
Now the Laplace operator
\begin{align*}
\Delta=\frac{\partial}{\partial x} \frac{\partial}{\partial
\bar{x}} - \sum_{i=1}^n \frac{\partial}{\partial \xi_i}
\frac{\partial}{\partial \bar{\xi}_i}: S^k(\C^{2|2n})\rightarrow
S^{k-2}(\C^{2|2n})
\end{align*}
is surjective of degree $-2$ for each $k\ge 0$. One checks that
$\Delta$ commutes with the action of ${\mathfrak g}$. This
establishes the following.

\begin{lemma} \label{mapdelta}
The map $\Delta: S^k(\C^{2|2n})\rightarrow S^{k-2}(\C^{2|2n})$ is
a surjective homomorphism of ${\mathfrak g}$-modules, and
$S^{k}(\C^{2|2n})/\ker\Delta\cong S^{k-2}(\C^{2|2n})$ as
${\mathfrak g}$-modules.
\end{lemma}

Consider the case $0\le k\le n$. In this case using the
combinatorial character formula of \cite[Theorem~3.7]{Lee} or
applying directly (\ref{BLformula}) we see that the character of
$\ker\Delta\subseteq S^{k}(\C^{2|2n})$ is equal to the character
of the irreducible module of highest weight $(k|0,\dots,0)$. This
immediately implies the following proposition.

\begin{proposition} \label{len}
For $0\le k\le n$ the ${\mathfrak g}$-module $\ker\Delta\subseteq
S^{k}(\C^{2|2n})$ is isomorphic to the irreducible highest weight
module of highest weight $(k|0,\dots,0)$.
\end{proposition}

Next we consider the case $k \ge 2n+1$. The following lemma is
easy to verify.

\begin{lemma}\label{dim>2n}
Let $k\ge 2n+1$.  Then ${\rm dim}S^{k}(\C^{2|2n})-{\rm
dim}S^{k-2}(\C^{2|2n})=2^{2n+1}$.
\end{lemma}

\begin{lemma}\label{hwv>2n}
Let $\Phi_{1}:=\sum_{i=1}^n\xi_i\bar{\xi}_i$ and $k\ge 2n$.  Set
\begin{align*}
\Gamma:=\sum_{i=0}^n(-1)^i{{k-n}\choose{i}}\bar{x}^{k-n-i}x^{n-i}\Phi^i_{1}
\in S^{k}(\C^{2|2n}).
\end{align*}
Then $\Gamma\not=0$ and we have $\Delta(\Gamma)=0$ and
$e_i\Gamma=0$, for $i=0,\cdots,n$.
\end{lemma}

\begin{proof}
Follows by a direct computation. \end{proof}

\begin{proposition} \label{ge2n+1}
Let $k\ge 2n+1$. There is an isomorphism of ${\mathfrak
g}$-modules:
$$\ker\Delta\cong L(k|0,\dots,0) \oplus L(2n-k|0,\dots,0).$$
\end{proposition}

\begin{proof}
It is clear that $x^k\in{\rm Ker\Delta}$ is a highest weight
vector of weight $(k|0,\dots, 0)$.  By Lemma~\ref{hwv>2n} the
irreducible $\fg$-module of highest weight $(2n-k|0,\cdots,0)$ is
also a composition factor of $\ker\Delta$. However, both weights
are typical, and hence the irreducible modules are equal to the
corresponding Kac modules which are of dimension $2^{2n}$. Now
Lemma \ref{dim>2n} implies that $\ker\Delta$ has only these two
composition factors. Finally, the weights $(k|0,\cdots,0)$ and
$(2n-k|0,\cdots,0)$ belong to different blocks and so indeed we
have a direct sum. \end{proof}

\begin{remark} \label{rem:autom}
Let $k\ge 2n+1$. Consider a new set of simple roots of $\fg$
associated with the following Dynkin diagram:
\begin{equation} %odd-dynkin%
\label{odd-dynkin} \setlength{\unitlength}{0.16in}
\begin{picture}(20,3.6)
\put(0,1){\makebox(0,0)[c]{$\bigcirc$}}
\put(3.4,1){\makebox(0,0)[c]{$\bigcirc$}}
\put(9.85,1){\makebox(0,0)[c]{$\bigcirc$}}
\put(13.35,-0.5){\makebox(0,0)[c]{$\bigotimes$}}
\put(13.35,2.5){\makebox(0,0)[c]{$\bigotimes$}}
\put(0.45,1){\line(1,0){2.45}} \put(3.82,1){\line(1,0){1.5}}
\put(5.7,1){$\ldots$} \put(7.2,1){\line(1,0){2.2}}
\put(10.28,1){\line(2,1){2.6}} \put(10.28,1){\line(2,-1){2.6}}
\put(13.28,2.05){\line(0,-1){2}}
\put(0,0){\makebox(0,0)[c]{$\delta_1-\delta_2$}}
\put(3.5,0){\makebox(0,0)[c]{$\delta_2-\delta_3$}}
\put(9.5,0){\makebox(0,0)[c]{$\delta_{n-1}-\delta_n$}}
\put(15.5,2.5){\makebox(0,0)[c]{$\delta_n-\epsilon$}}
\put(15.5,-0.5){\makebox(0,0)[c]{$\delta_n+\epsilon$}}
\end{picture}
%\end{array}$
\end{equation}
\vspace{.35cm}

\noindent Here $\delta_n \pm \epsilon$ are odd roots. We can show
via the method of odd reflections that the highest weights of the
two summands of $\ker \Delta$ in Proposition~\ref{ge2n+1} with
respect to this new Borel have Dynkin labels indicated as follows
(with the convention here and below that the unmarked ones are
$0$):

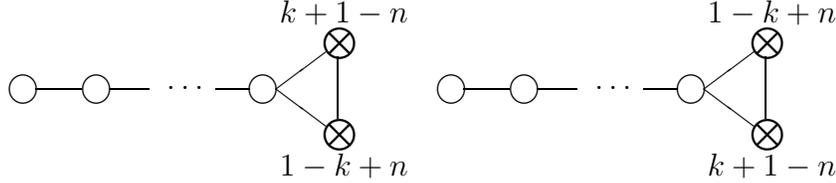
\begin{equation*}%k>2n%
%%\vspace*{-8ex}$
%\begin{array}{c c}
\setlength{\unitlength}{0.16in}
\begin{picture}(26,4)
\put(0,1){\makebox(0,0)[c]{$\bigcirc$}}
\put(2.4,1){\makebox(0,0)[c]{$\bigcirc$}}
\put(7.85,1){\makebox(0,0)[c]{$\bigcirc$}}
\put(10.35,-0.5){\makebox(0,0)[c]{$\bigotimes$}}
\put(10.35,2.5){\makebox(0,0)[c]{$\bigotimes$}}
\put(0.4,1){\line(1,0){1.5}} \put(2.82,1){\line(1,0){1.3}}
\put(4.7,1){$\ldots$} \put(6.3,1){\line(1,0){1.1}}
\put(8.28,1){\line(4,3){1.65}} \put(8.28,1){\line(4,-3){1.65}}
\put(10.28,2.05){\line(0,-1){2}}
\put(10.5,3.5){\makebox(0,0)[c]{$k+1-n$}}
\put(10.5,-1.5){\makebox(0,0)[c]{$1-k+n$}}
\put(14,1){\makebox(0,0)[c]{$\bigcirc$}}
\put(16.4,1){\makebox(0,0)[c]{$\bigcirc$}}
\put(21.85,1){\makebox(0,0)[c]{$\bigcirc$}}
\put(24.35,-0.5){\makebox(0,0)[c]{$\bigotimes$}}
\put(24.35,2.5){\makebox(0,0)[c]{$\bigotimes$}}
\put(14.4,1){\line(1,0){1.5}} \put(16.82,1){\line(1,0){1.3}}
\put(18.7,1){$\ldots$} \put(20.3,1){\line(1,0){1.1}}
\put(22.28,1){\line(4,3){1.65}} \put(22.28,1){\line(4,-3){1.65}}
\put(24.28,2.05){\line(0,-1){2}}
%
%\put(0,0){\makebox(0,0)[c]{$\beta_0$}}
%\put(2.4,0){\makebox(0,0)[c]{$\beta_1$}}
%\put(8.2,0){\makebox(0,0)[c]{$\beta_{n-2}$}}
%\put(12,2.5){\makebox(0,0)[c]{$\beta_{n-1}$}}
%\put(11.5,-0.5){\makebox(0,0)[c]{$\beta_{n}$}}
%
\put(24.5,3.5){\makebox(0,0)[c]{$1-k+n$}}
\put(24.5,-1.5){\makebox(0,0)[c]{$k+1-n$}}
%\put(0,4){\makebox(0,0){$k>2n$}}
\end{picture}
%\end{array}$
\end{equation*}
\vspace{.8cm}

\noindent Note that they are related via a Dynkin diagram
automorphism.
\end{remark}

It remains to consider the case $n+1\le k\le 2n$. Denote by $\e_l$
the irreducible character of $\text{sp}(2n)$ of highest weight
$\sum_{i=1}^l\delta_i$. The proofs of the following two lemmas are
straightforward and omitted.

\begin{lemma}\label{aux101}
Let $\la=\sum_{i=1}^nk_i\delta_i$, where $k_i=0,1$, for all $i$.
Suppose there exists $k_i=0$ and $k_j=1$ with $i<j$.  Let $x$ be
an indeterminate. Then we have
\begin{align*}
\sum_{w\in W}(-1)^{l(w)}w\left( e^{\la+\rho_0} \prod_{i=1}^n(1+ x
e^{-\delta_i})\right) =0.
\end{align*}
\end{lemma}

\begin{lemma}\label{aux102}
Let $\la=\sum_{i=1}^l\delta_i$, and $l\le n$. Let $x$ be an
indeterminate. Then we have
\begin{align*}
\frac{1}{D_0}\sum_{w\in W}(-1)^{l(w)}w\left( e^{\la+\rho_0}
\prod_{i=1}^n (1+ x e^{-\delta_i})\right) &=\sum_{j=0}^l x^{l -j}
\e_j.
\end{align*}
\end{lemma}

\begin{corollary}\label{aux103}
For $n+1\le k\le 2n$, the character ${\rm ch}
L(-\epsilon+\sum_{i=1}^{2n-k+1}\delta_i)$ equals
\begin{align}\label{aux106}
& \e_0\big{(}e^{-k\epsilon}+e^{(-k+2)\epsilon}+\cdots
e^{(-k+2(k-n-1))\epsilon}{)}\nonumber\\
&+\e_1\big{(}e^{(-k+1)\epsilon}+e^{(-k+3)\epsilon}+\cdots+
e^{(-k+1+2(k-n-1))\epsilon}{)} +\dots +\e_n e^{(n-k)\epsilon}.
\end{align}
\end{corollary}

\begin{proof}
By (\ref{BLformula}), the character ${\rm ch}
L(-\epsilon+\sum_{i=1}^{2n-k+1}\delta_i)$ is equal to
\begin{align*}
\frac{1}{D_0}&\sum_{w\in W}(-1)^{l(w)}w\Big{(}
e^{-\epsilon+\sum_{i=1}^{2n-k+1}\delta_i+\rho_0}\ \times\\
&(1+e^{-\epsilon+\delta_n})\cdots(1+e^{-\epsilon+\delta_{2n-k+2}})
\prod_{i=1}^n(1+e^{-\epsilon-\delta_i})\Big{)},
\end{align*}
which can then be written by Lemmas \ref{aux101} and \ref{aux102}
(with $x =e^{-\epsilon}$) as
\begin{align*}
%{\rm ch} L(-\epsilon+\sum_{i=1}^{2n-k+1}\delta_i)=
&e^{-\epsilon}\Big{(}
e^{(-2n+k-1)\epsilon}\e_0+e^{(-2n+k)\epsilon}\e_1+\cdots+\e_{2n-k+1}
\Big{)} \\
&+e^{-2\epsilon}\Big{(}
e^{(-2n+k-2)\epsilon}\e_0+e^{(-2n+k-1)\epsilon}\e_1+\cdots+\e_{2n-k+2}
\Big{)} +\dots\\
&+e^{(n-k)\epsilon}\Big{(}
e^{-n\epsilon}\e_0+e^{(-n+1)\epsilon}\e_1+\cdots+\e_{n} \Big{)}.
\end{align*}
The corollary now follows by collecting the coefficients of the
$\e_i$. \end{proof}

\begin{proposition} \label{le2n}
Let $n+1\le k\le 2n$.  We have
\begin{align*}
%{\rm ch}(S^k(\C^{2|2n}))-{\rm ch}S^{k-2}(\C^{2|2n})
{\rm ch} \ker \Delta ={\rm ch}L(k\epsilon) +{\rm
ch}L((2n-k)\epsilon) + {\rm
ch}L(-\epsilon+\sum_{i=1}^{2n-k+1}\delta_i).
\end{align*}
\end{proposition}

\begin{proof} Note that
$$
{\rm ch}K(k\epsilon)= e^{k\epsilon}\prod_{i=1}^n
(1+e^{-\epsilon+\delta_i})(1+e^{-\epsilon-\delta_i}),$$
and it can be rewritten as
\begin{align}\label{aux105}
%{\rm ch}K(k\epsilon)=
e^{k\epsilon}\Big{(}&
\e_0+e^{-\epsilon}\e_1+e^{-2\epsilon}(\e_2+\e_0)+e^{-3\epsilon}(\e_3+\e_1)+\cdots
+ e^{-n\epsilon}(\e_n+\e_{n-2}+\cdots)\nonumber\\
&+e^{-2n\epsilon}\e_0+ e^{(-2n+1)\epsilon}\e_1 +
e^{(-2n+2)\epsilon}(\e_2+\e_0)+e^{(-2n+3)\epsilon}(\e_3+\e_1)+\cdots\nonumber\\
&+\cdots+e^{-(n+1)\epsilon}(\e_{n-1}+\e_{n-3}+\cdots)\Big{)}.
\end{align}
Now a straightforward calculation shows that
\begin{align}\label{aux104}
%{\rm ch}(S^k(\C^{2|2n}))-{\rm ch}S^{k-2}(\C^{2|2n})
{\rm ch} \ker \Delta & =
(e^{(k-n)\epsilon}+e^{(n-k)\epsilon})\e_n\nonumber \\
&+ (e^{(k+1-n)\epsilon}+ e^{(k-1-n)\epsilon}+ e^{(n-k+1)\epsilon}
+e^{(n-k-1)\epsilon})\e_{n-1} + \cdots\nonumber\\
&+(e^{2(k-n-1)\epsilon}+\cdots+\widehat{1}
 +\cdots+e^{-2(k-n-1)\epsilon} )\e_{2n+2-k} +\cdots \\
&+(e^{(k-1)\epsilon}+e^{(k-3)\epsilon}+\cdots+e^{-(k-1)\epsilon})\e_1\nonumber\\
&+(e^{k\epsilon}+ e^{(k-2)\epsilon}+\cdots+e^{-k\epsilon})\e_0
\nonumber
\end{align}
where $\hat{1 }$ means as usual omission. One checks that
$(\ref{aux104})=(\ref{aux106})+(\ref{aux105})$. The proposition
now follows from the fact that ${\rm ch}K(k\epsilon)={\rm
ch}L(k\epsilon)+{\rm ch}L((2n-k)\epsilon)$. \end{proof}

Clearly, $S^0(\C^{2|2n})\cong \C$ and
$S^1(\C^{2|2n})\cong\C^{2|2n}$. The composition factors of
$S^k(\C^{2|2n})$ for every $k$ are now described explicitly by
combining Lemma~\ref{mapdelta}, Propositions~\ref{len},
\ref{ge2n+1} and \ref{le2n}.

\begin{remark}
Let $n+1\le k\le 2n$. We can show via the method of odd
reflections that the highest weights of the three summands of
$\ker \Delta$ in Proposition~\ref{le2n} with respect to the set of
simple roots (\ref{odd-dynkin})  have Dynkin labels as indicated
below:

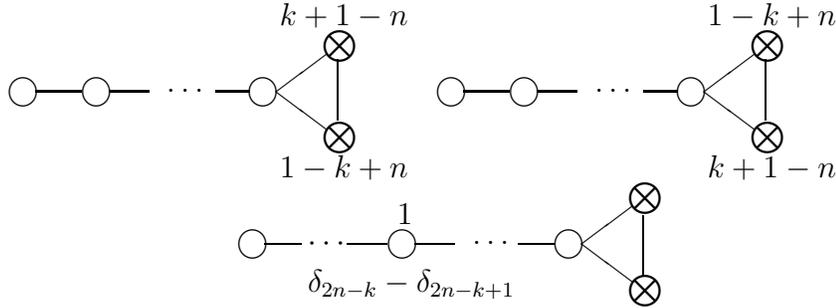
\begin{equation*}%n<k\le2n%
\setlength{\unitlength}{0.16in}
\begin{picture}(26,10)
\put(0,6){\makebox(0,0)[c]{$\bigcirc$}}
\put(2.4,6){\makebox(0,0)[c]{$\bigcirc$}}
\put(7.85,6){\makebox(0,0)[c]{$\bigcirc$}}
\put(10.35,4.5){\makebox(0,0)[c]{$\bigotimes$}}
\put(10.35,7.5){\makebox(0,0)[c]{$\bigotimes$}}
\put(0.4,6){\line(1,0){1.5}} \put(2.82,6){\line(1,0){1.3}}
\put(4.7,6){$\ldots$} \put(6.3,6){\line(1,0){1.1}}
\put(8.28,6){\line(4,3){1.65}} \put(8.28,6){\line(4,-3){1.65}}
\put(10.28,7.05){\line(0,-1){2}}
\put(10.5,8.5){\makebox(0,0)[c]{$k+1-n$}}
\put(10.5,3.5){\makebox(0,0)[c]{$1-k+n$}}
\put(14,6){\makebox(0,0)[c]{$\bigcirc$}}
\put(16.4,6){\makebox(0,0)[c]{$\bigcirc$}}
\put(21.85,6){\makebox(0,0)[c]{$\bigcirc$}}
\put(24.35,4.5){\makebox(0,0)[c]{$\bigotimes$}}
\put(24.35,7.5){\makebox(0,0)[c]{$\bigotimes$}}
\put(14.4,6){\line(1,0){1.5}} \put(16.82,6){\line(1,0){1.3}}
\put(18.7,6){$\ldots$} \put(20.3,6){\line(1,0){1.1}}
\put(22.28,6){\line(4,3){1.65}} \put(22.28,6){\line(4,-3){1.65}}
\put(24.28,7.05){\line(0,-1){2}}
%
%\put(0,0){\makebox(0,0)[c]{$\beta_0$}}
%\put(2.4,0){\makebox(0,0)[c]{$\beta_1$}}
%\put(8.2,0){\makebox(0,0)[c]{$\beta_{n-2}$}}
%\put(12,2.5){\makebox(0,0)[c]{$\beta_{n-1}$}}
%\put(11.5,-0.5){\makebox(0,0)[c]{$\beta_{n}$}}
%
\put(24.5,8.5){\makebox(0,0)[c]{$1-k+n$}}
\put(24.5,3.5){\makebox(0,0)[c]{$k+1-n$}}
%\put(0,3){\makebox(0,0){$n<k\le 2n$}}
%
\put(7.5,1){\makebox(0,0)[c]{$\bigcirc$}}
\put(10,1){\makebox(0,0)[c]{$\cdots$}}
\put(12.4,1){\makebox(0,0)[c]{$\bigcirc$}}
\put(17.85,1){\makebox(0,0)[c]{$\bigcirc$}}
\put(20.35,-0.5){\makebox(0,0)[c]{$\bigotimes$}}
\put(20.35,2.5){\makebox(0,0)[c]{$\bigotimes$}}
\put(7.9,1){\line(1,0){1.2}} \put(10.6,1){\line(1,0){1.3}}
\put(12.82,1){\line(1,0){1.3}} \put(14.7,1){$\ldots$}
\put(16.2,1){\line(1,0){1.15}} \put(18.28,1){\line(4,3){1.65}}
\put(18.28,1){\line(4,-3){1.65}} \put(20.28,1.95){\line(0,-1){2}}
\put(12.5,2){\makebox(0,0)[c]{$1$}}
\put(12.7,-0.3){\makebox(0,0)[c]{$\delta_{2n-k}-\delta_{2n-k+1}$}}
%\put(0,0){\makebox(0,0)[c]{$\beta_0$}}
%\put(2.4,0){\makebox(0,0)[c]{$\beta_1$}}
%\put(8.2,0){\makebox(0,0)[c]{$\beta_{n-2}$}}
%\put(12,2.5){\makebox(0,0)[c]{$\beta_{n-1}$}}
%\put(11.5,-0.5){\makebox(0,0)[c]{$\beta_{n}$}}
\end{picture}
\vspace{.6cm}
\end{equation*}
Note that all three weights are in the same block and two of them
are related by a diagram automorphism.
\end{remark}

For the sake of completeness, we remark that $\ker \Delta$ with $k
\le n$ (see Proposition~\ref{len}) with respect to the new Dynkin
diagram (\ref{odd-dynkin}) have the following Dynkin labels:

\vspace{0.5cm}
\begin{equation*}%k\le n%
%\vspace*{-8ex}$
%\begin{array}{c c}
\setlength{\unitlength}{0.16in}
\begin{picture}(20,3)
\put(15,1){\makebox(0,0)[c]{$\bigcirc$}}
\put(17.4,1){\makebox(0,0)[c]{$\bigcirc$}}
\put(22.85,1){\makebox(0,0)[c]{$\bigcirc$}}
\put(25.35,-0.5){\makebox(0,0)[c]{$\bigotimes$}}
\put(25.35,2.5){\makebox(0,0)[c]{$\bigotimes$}}
\put(15.45,1){\line(1,0){1.45}} \put(17.82,1){\line(1,0){1.3}}
\put(19.7,1){$\ldots$} \put(21.2,1){\line(1,0){1.15}}
\put(23.28,1){\line(4,3){1.65}} \put(23.28,1){\line(4,-3){1.65}}
\put(25.28,2.05){\line(0,-1){2}}
%
%\put(0,0){\makebox(0,0)[c]{$\beta_0$}}
%\put(2.4,0){\makebox(0,0)[c]{$\beta_1$}}
%\put(8.2,0){\makebox(0,0)[c]{$\beta_{n-2}$}}
%\put(12,2.5){\makebox(0,0)[c]{$\beta_{n-1}$}}
%\put(11.5,-0.5){\makebox(0,0)[c]{$\beta_{n}$}}
%
\put(25.4,3.6){\makebox(0,0)[c]{$1$}}
\put(25.4,-1.5){\makebox(0,0)[c]{$1$}}
\put(16,4){\makebox(0,0){$(k=n)$}}
\put(-3.5,1){\makebox(0,0)[c]{$\bigcirc$}}
\put(-1,1){\makebox(0,0)[c]{$\cdots$}}
\put(1.4,1){\makebox(0,0)[c]{$\bigcirc$}}
\put(6.85,1){\makebox(0,0)[c]{$\bigcirc$}}
\put(9.35,-0.5){\makebox(0,0)[c]{$\bigotimes$}}
\put(9.35,2.5){\makebox(0,0)[c]{$\bigotimes$}}
\put(-3.1,1){\line(1,0){1.2}} \put(-0.4,1){\line(1,0){1.3}}
\put(1.82,1){\line(1,0){1.3}} \put(3.7,1){$\ldots$}
\put(5.2,1){\line(1,0){1.15}} \put(7.28,1){\line(4,3){1.65}}
\put(7.28,1){\line(4,-3){1.65}} \put(9.28,1.95){\line(0,-1){2}}
\put(1.5,2){\makebox(0,0)[c]{$1$}}
\put(1.7,-0.3){\makebox(0,0)[c]{$\delta_{k}-\delta_{k+1}$}}
\put(-2,4){\makebox(0,0){$(k<n)$}}
%\put(0,0){\makebox(0,0)[c]{$\beta_0$}}
%\put(2.4,0){\makebox(0,0)[c]{$\beta_1$}}
%\put(8.2,0){\makebox(0,0)[c]{$\beta_{n-2}$}}
%\put(12,2.5){\makebox(0,0)[c]{$\beta_{n-1}$}}
%\put(11.5,-0.5){\makebox(0,0)[c]{$\beta_{n}$}}
\end{picture}
\vspace{.6cm}
\end{equation*}

\end{document}